\documentstyle[11pt,twoside]{article}
\include{graphicx}
\title{\bf An asymptotic expansion for the error term in the Brent-McMillan algorithm for Euler's constant}
\author{\sc R. B. Paris\footnote{E-mail address:\ \ {\tt r.paris@abertay.ac.uk}}\\
\\
{\em Division of Computing and Mathematics,}\\
{\em Abertay University, Dundee DD1 1HG, UK}\\
}

\begin{document}
\newcommand{\bee}{\begin{equation}}
\newcommand{\ee}{\end{equation}}
\def\f#1#2{\mbox{${\textstyle \frac{#1}{#2}}$}}
\def\dfrac#1#2{\displaystyle{\frac{#1}{#2}}}
\newcommand{\fr}{\frac{1}{2}}
\newcommand{\fs}{\f{1}{2}}
\newcommand{\g}{\Gamma}
\newcommand{\br}{\biggr}
\newcommand{\bl}{\biggl}
\newcommand{\ra}{\rightarrow}
\renewcommand{\topfraction}{0.9}
\renewcommand{\bottomfraction}{0.9}
\renewcommand{\textfraction}{0.05}
\newcommand{\mcol}{\multicolumn}
\date{}
\maketitle
\pagestyle{myheadings}
\markboth{\hfill {\it R.B. Paris} \hfill}
{\hfill {\it Error expansion for the Brent-McMillan algorithm} \hfill}
\begin{abstract} 
The Brent-McMillan algorithm is the fastest known procedure for the high-precision computation of Euler's constant $\gamma$ and is based on the modified Bessel functions $I_0(2x)$ and $K_0(2x)$. An error estimate for this algorithm relies on the optimally truncated asymptotic expansion for the product $I_0(2x) K_0(2x)$ when $x$ assumes large positive integer values. An asymptotic expansion for this optimal error term is derived by exploiting the techniques developed in hyperasymptotics, thereby enabling more precise information on the error term than recently obtained bounds and estimates.
\vspace{0.4cm}

\vspace{0.3cm}

\end{abstract}

\vspace{0.2cm}

\noindent $\,$\hrulefill $\,$

\vspace{0.2cm}

\begin{center}
{\bf 1. \  Introduction}
\end{center}
\setcounter{section}{1}
\setcounter{equation}{0}
\renewcommand{\theequation}{\arabic{section}.\arabic{equation}}
The Brent-McMillan algorithm \cite{BM} (when implemented with binary splitting) is the fastest known method of high-precision computation of Euler's constant $\gamma$. This relies on the formula \cite[(10.31.3)]{DLMF}
\bee\label{e11}
\gamma=\frac{S_0(2x)}{I_0(2x)}-\log\,x-\frac{K_0(2x)}{I_0(2x)},
\ee
where throughout we take $x$ to be a positive integer, $I_0(x)$, $K_0(x)$ are the standard modified Bessel functions and
\[S_0(x):=\sum_{n=0}^\infty \frac{H_n}{(n!)^2}\,(\fs x)^{2n},\qquad H_n:=1+\frac{1}{2}+\cdots +\frac{1}{n}.\]
For large $x$, the final term in (\ref{e11}) is $O(e^{-4x})$. Greater precision can be achieved following the suggestion made in \cite{BM} of truncating the asymptotic expansion
\bee\label{e12}
I_0(2x) K_0(2x)\sim\frac{1}{4x}\sum_{k=0}^\infty \frac{((2k)!)^3}{(k!)^4 (16x)^{2k}}\qquad (x\to+\infty)
\ee
at its optimal truncation index $k=2x$ (corresponding to truncation at, or near, the least term), followed by computing the term $K_0(2x)/I_0(2x)$ from $I_0(2x) K_0(2x)/(I_0(2x))^2$.
In \cite{BJ}, Brent and Johansson obtained a bound for the remainder term in the optimally truncated expansion (\ref{e12}) given by $24e^{-8x}$, thereby providing rigour to the algorithm. More recently, Demailly \cite{D} established the leading
large-$x$ behaviour of this remainder, together with an error bound, in the form
\bee\label{e12a}
-e^{-4x}\bl(\frac{5x^{-3/2}}{24\sqrt{2\pi}}+\epsilon(x)\br),\qquad |\epsilon(x)|<\frac{0.863}{x^2}.
\ee
This leads to the error in the optimally truncated expansion of the final term in (\ref{e11}) given by $-5\sqrt{2\pi}e^{-8x}/(12x^{1/2})$
to leading order.

The problem with using the well-known asymptotic expansions of the modified Bessel functions is that the positive real axis is a Stokes line for $I_0(x)$ (but not for $K_0(x)$). The standard expansion \cite[(10.40.5)]{DLMF}, \cite[p.~203]{W}
\bee\label{e13}
I_0(x)\sim \frac{e^x}{\sqrt{2\pi x}}\sum_{k=0}^\infty \frac{(\fs)_k (\fs)_k}{k! (2x)^k}+\frac{ie^{-x}}{\sqrt{2\pi x}} \sum_{k=0}^\infty \frac{(-)^k(\fs)_k (\fs)_k}{k! (2x)^k}\qquad (x\to+\infty)
\ee
is clearly inadequate, since this predicts a purely imaginary exponentially small contribution as $x\to+\infty$ when clearly it must be real. The correct form of the expansion of $I_0(x)$ for $x\to+\infty$ that takes into account the Stokes phenomenon on the positive $x$-axis has been considered in \cite{P17}. In this paper we derive an asymptotic expansion for the remainder in the optimally truncated expansion (\ref{e12}) by applying the first stage of the hyperasymptotic expansion process (also known as exponential improvement) to a suitable integral representation for the product $I_0(x) K_0(x)$, thereby bypassing the above-mentioned problem. A discussion of the new theory of hyperasymptotics (initiated by Berry \cite{Berry}) can be found in the book \cite[Ch. 6]{PK} in the context of the confluent hypergeometric functions; see also \cite[Section 2.11]{DLMF}. We present some numerical results to illustrate the accuracy of the expansion so obtained.

\vspace{0.6cm}

\begin{center}
{\bf 2. \ Exponentially improved expansion for $I_0(x) K_0(x)$}
\end{center}
\setcounter{section}{2}
\setcounter{equation}{0}
\renewcommand{\theequation}{\arabic{section}.\arabic{equation}}
We start with the Mellin-Barnes integral representation \cite[p.~116]{PK}
\[I_0(x)K_0(x)=\frac{1}{2\pi i}\int_{-c-\infty i}^{-c+\infty i} \g^3(s+\fs) \g(-s)\frac{\cos \pi s}{2\pi^{3/2}} x^{-2s-1}ds\qquad (|\arg\,x|<\fs\pi,\ 0<c<\fs).\]
Throughout this paper we shall restrict $x$ to be a positive integer in keeping with the strategy of the Brent-McMillan algorithm, although the analysis can be developed for complex $x$. The integrand has simple poles situated at $s=0, 1, 2, \ldots$ and double poles at $s=-\fs, -\f{3}{2}, \ldots\,$. 

We consider the integral taken round the rectangular contour with vertices at $c\pm iT$ and $N-c'\pm iT$, where $N$ is (for the moment) an arbitrary positive integer and $0<c'<1$. Use of the well-known approximation (with $\sigma$ real) $\g(\sigma\pm it)=O(t^{\sigma-1/2} e^{-\frac{1}{2}\pi t})$ as $t\to+\infty$, shows that the contribution from the upper and lower sides $s=\sigma\pm iT$, $c\leq\sigma\leq N-c'$ vanishes as $T\to\infty$, since the modulus of the integrand on these paths is $O(T^{2\sigma-1/2}x^{-2\sigma-1} e^{-\pi T})$.
Displacement of the integration path to the right over the first $N$ poles, together with the fact that the residue of $\g(-s)$ at $s=k$ is $(-1)^{k-1}/k!$, then shows that
\bee\label{e21}
I_0(x)K_0(x)=\frac{1}{2x}\sum_{k=0}^{N-1} \frac{((2k)!)^3}{(k!)^4 (8x)^{2k}}+R_N(x),
\ee
where the remainder $R_N(x)$ is 
\bee\label{e22}
R_N(x)=-\frac{1}{2\pi i}\int_{L_N}\frac{\g^3(s+\fs)}{\g(s+1)}\,\cot \pi s\,\frac{x^{-2s-1}}{2\sqrt{\pi}}ds,
\ee
and $L_N$ denotes the rectilinear path $(-c'\!+\!N\!-\!\infty i, -c'\!+\!N\!+\!\infty i)$.

We now choose $N$ to be the optimal truncation index of the expansion in (\ref{e21}), which is easily verified to be $N=x$. As a consequence, since $x\to+\infty$ the variable $s$ in the quotient of gamma functions in (\ref{e22}) is uniformly large on the displaced path $L_N$. From Lemma 2.2 in \cite[p.~39]{PK} we have the inverse factorial expansion
\bee\label{e22b}
\frac{\g^3(s+\fs)}{\g(s+1)}=\frac{2\sqrt{\pi}}{2^{2s}}\bl\{\sum_{j=0}^{M-1} (-)^j c_j \g(2s-j)+\rho_M(s) \g(2s-M)\br\}
\ee
for positive integer $M$, where
\bee\label{e22a}
c_0=1,\quad c_1=\fs,\quad c_2=\f{5}{8},\quad c_3=\f{21}{16},\quad c_4=\f{507}{128},\quad c_5=\f{4035}{256};
\ee
see the appendix. The remainder function $\rho_M(s)$ is analytic in $s$ except at the points $s=-\fs, -\f{3}{2}, \ldots$ and is such that $\rho_M(s)=O(1)$ for $|s|\to\infty$ in $|\arg\,s|<\pi$. Then we obtain
\[R_N(x)=-\frac{1}{x}\bl\{\sum_{j=0}^{M-1}(-)^jc_j\,\frac{1}{2\pi i}\int_{L_N}\g(2s-j) \cot \pi s \,(2x)^{-2s}ds+R_{M,N}(x)\br\},\]
where
\begin{eqnarray*}
R_{M,N}(x)&=&\frac{1}{2\pi i}\int_{L_N} \rho_M(s) \g(2s-M) \cot \pi s\,(2x)^{-2s} ds\\
&=&\frac{1}{2\pi i}\int_{L_N}\rho_M(s) \frac{ \g(2s-M)}{\sin 2\pi s}\,(1+\cos 2\pi s)\,(2x)^{-2s}ds.
\end{eqnarray*}

The remainder $R_{M,N}(x)$ can be split into three separate integrals with variables $(2x)^{-2s}$ and $(2xe^{\pm \pi i})^{-2s}$, to each of which we can apply Lemma 2.9 in \cite[p.~75]{PK}. Since we have chosen $N=x$, we therefore obtain the order estimates $O(x^{-M-\frac{1}{2}}e^{-2x})$ for the integral with variable $(2x)^{-2s}$ and $O(x^{-M}e^{-2x})$ for the integrals with variables $(2xe^{\pm\pi i})^{-2s}$. Hence
$R_{M,N}(x)=O(x^{-M}e^{-2x})$ as $x\to+\infty$. Then we find
\bee\label{e23}
R_N(x)=-\frac{1}{x}\bl\{\sum_{j=0}^{M-1}(-)^j c_j \frac{1}{2\pi i}\int_{L_N} \frac{\g(2s-j)}{\sin 2\pi s}\,(1+\cos 2\pi s)\,(2x)^{-2s}ds+O(x^{-M}e^{-2x})\br\}.
\ee

We now introduce the so-called {\it terminant function} $T_\nu(z)$ defined\footnote{In \cite[(2.11.11)]{DLMF} this function is denoted by $F_\nu(z)$ and is expressed as a multiple of the exponential integral $E_\nu(z)=z^{\nu-1}\g(1-\nu,z)$.} as a multiple of the incomplete gamma function $\g(a,z)$ by
\[{T}_\nu(z):=\frac{\g(\nu)}{2\pi}\, \g(1-\nu,z).\]
 From the formula connecting $\g(a,ze^{\pm\pi i})$ given in \cite[(8.2.10)]{DLMF}
we have the connection formula (compare also \cite[(6.2.45)]{PK})
\bee\label{e24}
T_\nu(ze^{-\pi i})=e^{2\pi i\nu}\{T_\nu(ze^{\pi i})-ie^{-\pi i\nu}\}.
\ee
The Mellin-Barnes integral representation of this function is \cite[(6.2.7)]{PK}
\bee\label{e25}
-2z^\nu e^z T_\nu(z)=\frac{1}{2\pi i}\int_{-c-\infty i}^{-c+\infty i} \frac{\g(s+\nu)}{\sin \pi s}\,z^{-s}ds\qquad (|\arg\,z|<\f{3}{2}\pi,\ 0<c<1)
\ee
provided $\nu\neq 0, -1, -2, \ldots\,$. Then, if we make the change of variable $s\to s+N$ in the integrals appearing in (\ref{e23}), write $\cos 2\pi s$ in terms of exponentials, and use (\ref{e25}) (when it is supposed that $M<2N$) these integrals can be written as
\[\frac{1}{2\pi i} \int_{-c-\infty i}^{-c+\infty i} \frac{\g(2s+2N-j)}{\sin 2\pi s}\,(1+\cos 2\pi s)\,(2x)^{-2s-2N}ds\hspace{4cm}\]
\[\hspace{2.2cm}=-(2x)^{-j}\bl\{e^{2x}T_{2N-j}(2x)+\fs e^{-2x} (-)^j[T_{2N-j}(2xe^{\pi i})+T_{2N-j}(2xe^{-\pi i})]\br\}\]
\[=-(2x)^{-j}\bl\{e^{2x}T_{2N-j}(2x)+e^{-2x}[(-)^j T_{2N-j}(2xe^{\pi i})-\fs i]\br\}\]
upon application of (\ref{e24}).

This then yields the expansion for $R_N(x)$ given by
\bee\label{e26}
R_N(x)=\frac{e^{-2x}}{x}\bl\{\sum_{j=0}^{M-1}\frac{(-)^j c_j}{(2x)^j}\bl\{e^{4x}T_{2N-j}(2x)+(-)^j T_{2N-j}(2xe^{\pi i})-\fs i\br\}+O(x^{-M})\br\}.
\ee
It now remains to exploit the known asymptotic expansions of the terminant function $T_\nu(x)$ when $\nu\sim x$ as $x\to+\infty$, which we carry out in the next section.

\vspace{0.6cm}

\begin{center}
{\bf 3. \  An asymptotic expansion for $R_N(x)$}
\end{center}
\setcounter{section}{3}
\setcounter{equation}{0}
\renewcommand{\theequation}{\arabic{section}.\arabic{equation}}
The asymptotic expansion of the terminant function $T_\nu(z)$ for large $\nu$ and complex $z$, when $\nu\sim |z|$, has been discussed in detail by Olver in \cite{O91}; see also \cite[Section 2.11]{DLMF} and the detailed account in \cite[pp.~259--265]{PK}. By expressing $T_\nu(z)$ in terms of the Laplace integral
\[T_\nu(z)=\frac{e^{-z}}{2\pi}\int_0^\infty e^{-zt}\,\frac{t^{\nu-1}}{1+t}\, dt,\]
Olver established by application of the saddle-point method that when $\mu\sim x$ (and bounded integer $j$)
\bee\label{e31}
T_{\mu-j}(x)=\frac{e^{-2x}}{2\sqrt{2\pi x}}\bl\{\sum_{k=0}^{K-1}A_{k,j}x^{-k}+O(x^{-K})\br\}\qquad (x\to+\infty),
\ee
where $A_{0,j}=1$ ($j\geq 0$) and
\begin{eqnarray}
A_{1,j}&=&\f{1}{6}(2-6\gamma_j+3\gamma_j^2),\quad A_{2,j}=\f{1}{288}(-11-120\gamma_j+300\gamma_j^2-192\gamma_j^3+36\gamma_j^4),\nonumber\\
A_{3,j}&=&\f{2}{51840}(-587+3510\gamma_j+9765\gamma_j^2-26280\gamma_j^3+18900\gamma_j^4-5400\gamma_j^5+540\gamma_j^6),\nonumber\\
A_{4,j}&=&\f{1}{2448320}(120341-44592\gamma_j-521736\gamma_j^2-722880\gamma_j^3+2336040\gamma_j^4-1826496\gamma_j^5\nonumber\\
&&\hspace{5cm}+635040\gamma_j^6-103680\gamma_j^7+6480\gamma_j^8),\label{e31a}
\end{eqnarray}
with 
\bee\label{e31b}
\gamma_j:=\mu-x-j\qquad (0\leq j\leq K-1).
\ee

On the negative real axis, where a saddle point and a simple pole become coincident in the above Laplace integral, we have the expansion
\bee\label{e32}
T_{\mu-j}(xe^{\pi i})=e^{-\pi i\mu}(-)^j\bl\{\fs i+\frac{1}{\sqrt{2\pi x}}\bl(\sum_{k=0}^{K-1} (\fs)_k G_{2k,j}\,(\fs x)^{-k}+O(x^{-K})\br)\br\}\qquad (x\to+\infty),
\ee
where the coefficients $G_{k,j}$ result from the expansion
\[\frac{\tau^{\gamma_j-1}}{1-\tau}\,\frac{d\tau}{dw}=-\frac{1}{w}+\sum_{k=0}^\infty G_{k,j}w^k,\qquad \fs w^2=\tau-\log\,\tau-1.\]
The branch of $w(\tau)$ is chosen such that $w\sim \tau-1$ as $\tau\ra 1$. Upon reversion of the $w$-$\tau$ mapping to yield
\[\tau=1+w+\f{1}{3}w^2+\f{1}{36}w^3-\f{1}{270}w^4+\f{1}{4320}w^5+ \cdots\ ,\]
it is found with the help of {\it Mathematica} that the first five even-order coefficients $G_{2k,j}\equiv 6^{-2k} {\hat G}_{2k,j}$ are
\begin{eqnarray}
{\hat G}_{0,j}\!\!&=&\!\!\f{2}{3}-\gamma_j,\qquad {\hat G}_{2,j}=\f{1}{15}(46-225\gamma_j+270\gamma_j^2-90\gamma_j^3), \nonumber\\
{\hat G}_{4,j}\!\!&=&\!\!\f{1}{70}(230-3969\gamma_j+11340\gamma_j^2-11760\gamma_j^3+5040\gamma_j^4
-756\gamma_j^5),\nonumber\\
{\hat G}_{6,j}\!\!&=&\!\!\f{1}{350}(-3626-17781\gamma_j+183330\gamma_j^2-397530\gamma_j^3+370440\gamma_j^4
-170100\gamma_j^5\nonumber\\
&&\hspace{7cm}+37800\gamma_j^6-3240\gamma_j^7),\nonumber\\
{\hat G}_{8,j}\!\!&=&\!\!\f{1}{231000}(-4032746+43924815\gamma_j+88280280\gamma_j^2-743046480\gamma_j^3\nonumber\\
&&+1353607200\gamma_j^4-1160830440\gamma_j^5+541870560\gamma_j^6
-141134400\gamma_j^7\nonumber\\
&&\hspace{6cm}+19245600\gamma_j^8-1069200\gamma_j^9).\label{e32a}
\end{eqnarray}

Substitution of the expansions (\ref{e31}) and (\ref{e32}) with $\mu=2N$ into $R_N(x)$ in (\ref{e26}) then yields
(where we put $K=M$ for convenience)
\[R_N(x)=\frac{e^{-2x}}{4\sqrt{\pi}x^{3/2}}\sum_{j=0}^{M-1}\frac{(-)^jc_j}{(2x)^j} \sum_{k=0}^{M-1}\frac{D_{k,j}}{(2x)^k}+O(x^{-M-1}e^{-2x})\]
for $x\to+\infty$, where
\bee\label{e33a}
D_{k,j}:=A_{k,j}+2^{k+1} (\fs)_k G_{2k,j}
\ee
and, since the variables in the terminant functions in (\ref{e26}) involve $2x=2N$, we have from (\ref{e31b}) that $\gamma_j=-j$. 
Then we obtain the following theorem:
\newtheorem{theorem}{Theorem}
\begin{theorem}$\!\!\!.$\  Let $M$, $N$ be positive integers and the variable $x$ assume integer values. Then the remainder $R_N(x)$ in the optimally truncated asymptotic expansion for $I_0(x) K_0(x)$ in (\ref{e21}) when $N=x$ has the expansion
\bee\label{e33}
R_N(x)=\frac{e^{-2x}}{4\sqrt{\pi} x^{3/2}}\bl\{\sum_{j=0}^{M-1}B_j (2x)^{-j}+O(x^{-M+\frac{1}{2}})\br\}
\ee
as $x\to+\infty$. The coefficients $B_j$ are defined by
\bee\label{e34}
B_j=\sum_{k=0}^j(-)^k c_k D_{j-k,k},
\ee
where the coefficients $c_k$ and $D_{k,j}$ are specified in (\ref{e22a}) and (\ref{e33a}). The quantities $A_{k,j}$ and $G_{2k,j}$ appearing in (\ref{e33a}) are defined in (\ref{e31a}) and (\ref{e32a}) with $\gamma_j=-j$.
\end{theorem}

Routine computations show that
\[B_0=\frac{7}{3},\qquad B_1=-\frac{449}{270},\qquad B_2=\frac{55949}{3024},\qquad B_3=-\frac{87499}{17010},\]
\[B_4=\frac{137885143760267}{7067908108800}.\]
This produces the expansion
\bee\label{e35}
R_N(x)\sim \frac{7e^{-2x}}{12\sqrt{\pi} x^{3/2}}\bl\{1-\frac{449}{1260\,x}+\frac{55949}{282240\,x^2}-\frac{87499}{317520\,x^3}+\frac{137885143760267}{263868569395200\,x^4}+\cdots \br\}
\ee
as $x\to+\infty$, which is the main result of the paper. In Table 1 we present values of the absolute relative error in the computation of the expansion for $R_N(x)$ in (\ref{e33}) for different $x$ and truncation index $M$ compared with the exact evaluation from (\ref{e21}).
\begin{table}[th]
\caption{\footnotesize{Values of the absolute relative error in the computation of $R_N(x)$ from (\ref{e33}).}} \label{t1}
\begin{center}
\begin{tabular}{|l||l|l|l|}
\hline
&&&\\[-0.3cm]
\mcol{1}{|c||}{$M$} & \mcol{1}{c|}{$x=50$} & \mcol{1}{c|}{$x=100$} & \mcol{1}{c|}{$x=150$}\\
\hline
&&&\\[-0.3cm]
1 & $7.100\times 10^{-3}$ & $3.557\times 10^{-3}$ & $2.373\times 10^{-3}$  \\
2 & $7.772\times 10^{-5}$ & $1.962\times 10^{-5}$ & $8.750\times 10^{-6}$  \\
3 & $2.140\times 10^{-6}$ & $2.714\times 10^{-7}$ & $8.082\times 10^{-8}$  \\
4 & $8.065\times 10^{-8}$ & $5.130\times 10^{-9}$ & $1.020\times 10^{-9}$  \\
5 & $3.555\times 10^{-9}$ & $1.137\times 10^{-10}$& $1.510\times 10^{-11}$  \\
[.1cm]\hline\end{tabular}
\end{center}
\end{table}
\vspace{0.6cm}

\begin{center}
{\bf 4. \  Concluding remarks}
\end{center}
\setcounter{section}{4}
\setcounter{equation}{0}
\renewcommand{\theequation}{\arabic{section}.\arabic{equation}}
In the Brent-McMillan algorithm we have, with $N=x$,
\[I_0(2x) K_0(2x)=\frac{1}{4x}\sum_{k=0}^{2N-1} \frac{((2k)!)^3}{(k!)^4 (16x)^{2k}}+R_{2N}(2x),\]
where from (\ref{e35}) 
\bee\label{e41}
R_{2N}(2x)\sim \frac{7e^{-4x}}{24\sqrt{2\pi} x^{3/2}}\bl\{1-\frac{449}{2520\,x}+ \frac{55949}{1128960\,x^2}-\frac{87499}{2540160\,x^3}+\frac{137885143760267}{4221897110323200\,x^4}+\cdots\br\}
\ee
for $x\to+\infty$. 
From (\ref{e13}), the expansion of $(I_0(2x))^2$ (upon neglecting the exponentially small contribution) is
\[(I_0(2x))^2\sim \frac{e^{4x}}{4\pi x}\bl\{1+\frac{1}{8x}+\frac{5}{128x^2}+\frac{21}{1024x^3}+\frac{507}{32768x^4}+\cdots\br\}\qquad(x\to+\infty).\]
Then, from (\ref{e41}), we obtain
\begin{theorem}$\!\!\!.$\  
The error resulting from $K_0(2x)/I_0(2x)$ in the Brent-McMillan algorithm in (\ref{e11}) at optimal truncation  has the expansion
\bee\label{e42a}
\frac{R_{2N}(2x)}{(I_0(2x))^2}=\frac{7\sqrt{2\pi} e^{-8x}}{12 x^{1/2}}\,\bl\{1-\frac{191}{630x}+\frac{18211}{376320\,x^2}-\frac{799201}{16257024\,x^3}+\frac{116774621369177}{4221897110323200\,x^4}+\cdots\br\}
\ee
as $x\to+\infty$.
\end{theorem}

In \cite{D}, Demailly defined his remainder function $\Delta(x)$ using the optimal truncation index $k=2N$, instead of $k=2N-1$, and wrote
\[\Delta(x):=I_0(2x) K_0(2x)-\frac{1}{4x}\sum_{k=0}^{2N} \frac{((2k)!)^3}{(k!)^4 (16x)^{2k}}\qquad (N=x).\]
The connection with our $R_{2N}(2x)$ is consequently given by
\bee\label{e42}
\Delta(x)=R_{2N}(2x)-\frac{1}{4x}\,\frac{((4x)!)^3}{((2x)!)^4 (16x)^{4x}}.
\ee
The expansion of the second term in (\ref{e42}) can be obtained by application of Stirling's formula (see (\ref{a3}) with $s=2x$) to yield
\[
\Delta(x)\sim -\frac{5e^{-4x}}{24\sqrt{2\pi} x^{3/2}}\bl\{1-\frac{1}{1800\,x}-\frac{45449}{806400\,x^2}+\frac{294911}{5806080\,x^3}+\cdots\br\}
\]
as $x\to+\infty$. It is clear that the error estimate in (\ref{e12a}) considerably overestimates the first-order correction to $\Delta(x)$ in the limit $x\to+\infty$.

\vspace{0.6cm}

\begin{center}
{\bf Appendix:  Determination of the coefficients $c_j$ in the expansion (\ref{e22b})}
\end{center}
\setcounter{section}{1}
\setcounter{equation}{0}
\renewcommand{\theequation}{\Alph{section}.\arabic{equation}}
Use of the duplication formula $\g(2s)=\pi^{-1/2}2^{2s-1}\g(s) \g(s+\fs)$ shows that the inverse factorial expansion (\ref{e22b}) can be written as
\bee\label{a1}
\frac{2^{2s}}{2\sqrt{\pi}}\,\frac{\g^3(s+\fs)}{\g(s+1)\g(2s)}=\frac{1}{s}\bl(\frac{\g(s+\fs)}{\g(s)}\br)^{\!\!2}=\sum_{j=0}^{M-1}\frac{c_j}{(1-2s)_j}+\frac{(-)^M\rho_M(s)}{(1-2s)_M}.
\ee
From \cite[(5.11.13)]{DLMF} (see also \cite[(2.2.32)]{PK}) we have
\[\frac{\g(s+\fs)}{\g(s)}=s^{1/2}\bl\{1-\frac{1}{8s}+\frac{1}{128s^2}+\frac{5}{1024s^3}-\frac{21}{32768s^4}-\frac{399}{262144s^5}+\cdots\br\}\qquad (s\to+\infty),\]
whence
\bee\label{a2}
\frac{1}{s}\bl(\frac{\g(s+\fs)}{\g(s)}\br)^{\!\!2}=1-\frac{1}{4s}+\frac{1}{32s^2}+\frac{1}{128s^3}-\frac{5}{2048s^4}-\frac{23}{8192s^5}+\cdots
\ee
\[=c_0-\frac{c_1}{2s}+(c_2-c_1)\frac{1}{4s^2}+(-c_1+3c_2-c_3)\frac{1}{8s^3}+(-c_1+7c_2-6c_3+c_4)\frac{1}{16s^4}\]
\[\hspace{8cm}+(-c_1+15c_2-25c_3+10c_4-c_5)\frac{1}{32s^5}+ \cdots\]
upon expansion of the right-hand side of (\ref{a1}) in inverse powers of $s$.
Comparison of the coefficients of corresponding powers of $s$ then yields the values
\[c_0=1,\quad c_1=\fs,\quad c_2=\f{5}{8},\quad c_3=\f{21}{16},\quad c_4=\f{507}{128},\quad c_5=\f{4035}{256}.\]

We can make use of (\ref{a2}) and the duplication formula to obtain the expansion of the quantity required in Section 4
\[\frac{1}{2s}\,\frac{((2s)!)^3}{(s!)^4 (8s)^{2s}}=\frac{\g^3(s+\fs)}{2\pi^{3/2} \g(s+1) s^{2s+1}}=\frac{\g(2s)}{\pi s(2s)^{2s}}\,\frac{1}{s}\bl(\frac{\g(s+\fs)}{\g(s)}\br)^{\!\!2}
\]
\[=\frac{\g(2s)}{\pi s (2s)^{2s}}\,\bl\{1-\frac{1}{4s}+\frac{1}{32s^2}+\frac{1}{128s^3}-\frac{5}{2048s^4}-\frac{23}{8192s^5}+\cdots\br\}.\]
Use of the well-known expansion \cite[(5.11.3)]{DLMF}
\[\g(z)=\sqrt{2\pi} z^{z-\frac{1}{2}}e^{-z} \bl\{1+\frac{1}{12z}+\frac{1}{288z^2}-\frac{139}{51840z^3}+\cdots\br\}\qquad (z\to+\infty)\]
then shows that
\bee\label{a3}
\frac{1}{2s}\,\frac{((2s)!)^3}{(s!)^4 (8s)^{2s}}=\frac{e^{-2s}}{\sqrt{\pi} s^{3/2}}\bl\{1-\frac{5}{24s}+\frac{25}{1152s^2}+\frac{3551}{414720s^3}+O(s^{-4})\br\}
\ee
as $s\to+\infty$.

\end{document}